\documentclass[a4paper]{article}

\usepackage[english]{babel}
\usepackage[utf8x]{inputenc}
\usepackage[T1]{fontenc}
\usepackage{lscape}
\usepackage[a4paper,top=3cm,bottom=2cm,left=3cm,right=3cm,marginparwidth=1.75cm]{geometry}

\usepackage{amsmath}
\usepackage{graphicx}
\usepackage[colorinlistoftodos]{todonotes}
\usepackage[colorlinks=true, allcolors=blue]{hyperref}
\title{A note on rank constrained solutions to linear matrix equations}
\author{Shravan Mohan}

\begin{document}
\maketitle

\begin{abstract}
This preliminary note presents a heuristic for determining rank constrained solutions to linear matrix equations (LME). The method proposed here is based on minimizing a non-convex quadratic functional, which will hence-forth be termed as the \textit{Low-Rank-Functional} (LRF). Although this method lacks a formal proof/comprehensive analysis, for example in terms of a probabilistic guarantee for converging to a solution, the proposed idea is intuitive and has been seen to perform well in simulations. To that end, many numerical examples are provided to corroborate the idea.   
\end{abstract}

\section{Introduction}
Consider the feasibility problem:
\begin{eqnarray}
\mbox{Find}~~ &X \in \mathcal{R}^{N\times N}\nonumber \\
\mbox{such~that}~~ &\mathcal{A}(X) = b, ~ X \succeq 0 \mbox{~~and~~} \mbox{Rank} (X) = k,
\label{main_prob}
\end{eqnarray}
where $\mathcal{A}$ is a linear transformation from $\mathcal{R}^{N\times N}$ to $\mathcal{R}^{m}$ and $b\in \mathcal{R}^{m\times 1}$. It is also assumed that $m < N^2$. The case where the rank of the solution is constrained to be 1 is of particular interest to many applications such as combinatorics and signal processing. In the following discussions, (i) $[e;~f]$ would mean the concatenated vector by appending the vector $f$ below $e$, (ii) $\mbox{diag}(X)$ would mean the column vector of the diagonal entries of the matrix $X$, (iii) $X(m:n,p:q)$ would the sub-matrix of $X$ comprised of the elements with row index between $m$ and $n$  and column index between $p$ and $q$, (iv) $X(:,i)$ would be the $i^{\rm th}$ column of $X$, $X(i,:)$ would be the $i^{\rm th}$ row of $X$, (v) $X(m:n,i)$ would be the $i^{\rm th}$ column with elements from row $m$ to row $n$ (similar interpretation would hold for $X(i,m:n)$), (vi) $X(:)$ or $\mbox{\textbf{vec}}(X)$  would mean the vector obtained by vectorizing $X$, and (vii) $\mbox{Sym}^N$ and $S_+^N$ are the spaces of $N\times N$ symmetric matrices and positive semidefinite matrices, respectively.  Continuing, a typical decision problem in combinatorics would have the following form:
\begin{eqnarray}
\mbox{Find}~~ &x \in \{0,1\}^{N\times 1} \nonumber\\
\mbox{such~that}~~ &\mathcal{C}x = b,
\label{combi_prob}
\end{eqnarray}
where $C\in \mathcal{R}^{m\times N}$ and $b\in \mathcal{R}^{m\times 1}$. A standard method to convert this into a form given in \eqref{main_prob} is to firstly recognize the simple fact that $z \in \{0,1\}$ if and only if $z(z-1) = 0$. Now suppose $X = [1; ~x][1; ~x]^{\top}$. By construction, the matrix $X$ is positive semidefinite and has unity rank. Moreover, the constraint $x_i(x_i-1) = 0$, where $x_i$ is the $i^{\mbox th}$ element of the vector  $x$, can be equivalently written as $\mbox{diag}(X) = X(:,1)$. With all these facts in place, the feasibility problem \eqref{combi_prob} can be written as
\begin{eqnarray}
\mbox{Find}~~ &X \in \mathcal{R}^{N\times N}\nonumber \\
\mbox{such~that}~~ &\mathcal{C}X_{(\small{2:N,1})} = b, ~\mbox{diag}(X) = X_{(:,1)}, ~X \succeq 0 \mbox{~~and~~} \mbox{Rank} (X) = 1.
\label{combi_prob_conversion}
\end{eqnarray}
A problem from the field of non-linear equations which can be converted to rank constrained feasibility problem is the following: 
\begin{eqnarray}
\mbox{Find}~~ &[\theta_1,~\theta_2,\cdots,~\theta_N] \in [-\pi,\pi]^{N\times 1}\nonumber \\
\mbox{such~that}~~ & Ax = b, \mbox{~~where~} x = \left[e^{\mbox{j}\theta_1},~e^{\mbox{j}\theta_2},\cdots,~e^{\mbox{j}\theta_N}\right]^\top,
\label{nl_prob1}
\end{eqnarray}
and $A\in \mathcal{R}^{m\times N}$ and $b\in \mathcal{R}^{m\times 1}$. Here it is assumed that the given system of equations has a solution. As in the example in combinatorics, suppose that $X = [1; ~x][1; ~x]^{\mbox H}$. Again by construction, $X \succeq 0$ and $X$ has unity rank. A minor difference, as compared to the previous example, is the fact that the $X$ here is defined on the complex field. Note that $\mbox{diag}(X)=\textbf{1}$, a column of ones. With this, the problem in \eqref{nl_prob1} can be cast as:
\begin{eqnarray}
\mbox{Find}~~ &X \in \mathcal{R}^{N\times N}\nonumber \\
\mbox{such~that}~~ &\mathcal{A}X(\small{2:N,1}) = b, ~\mbox{diag}(X) = \textbf{1}, ~X \succeq 0 \mbox{~~and~~} \mbox{Rank} (X) = 1.
\label{combi_prob_conversion}
\end{eqnarray}
It is also worth mentioning an instance where a non-unity rank comes of use is in a very elegant formulation of the Optimal Power Flow problem proposed in \cite{lavaei2012zero}. It is shown in \cite{lavaei2012zero} that a sufficient condition for a solution to the dual problem to be optimal to the primal is that a specific positive semidefinite matrix (affine in the optimization variables) has a null space of dimension 2. The aforementioned examples highlight the advantage of using same rank-constrained feasibility formulation to represent a 0/1 programming problem on one hand, and a continuous domain feasibility problem on the other. Many more applications can be found in the paper \cite{fazel2004rank,recht2010guaranteed}, which serve as the motivation for studying this problem in detail. In this paper, the general rank constrained problem is also considered, which is given by:
\begin{eqnarray}
\mbox{Find}~~ &X \in \mathcal{R}^{N\times M}\nonumber \\
\mbox{such~that}~~ &\mathcal{A}(X) = b \mbox{~~and~~} \mbox{Rank} (X) = k.
\label{main_prob_gen}
\end{eqnarray}
A related problem is the \textit{Minimum-Rank-Problem} (MRP) given by:
\begin{eqnarray}
\mbox{Find}~~ &X \in \mathcal{R}^{M\times N} \mbox{~~with ~the ~minimum~ rank} \nonumber \\
\mbox{such~that}~~ &\mathcal{A}(X) = b \mbox{~and~} X \succeq 0.
\label{min_rank_prob}
\end{eqnarray}
A popular heuristic employed to solve \eqref{min_rank_prob} is called the \textbf{log-det} heuristic, first proposed in the seminal paper \cite{fazel2003log}. This method relies on the intuition that minimizing $\log(\det(X))$ would naturally reduce the singular values and hence lead to rank minimization. In particular, this is achieved by the following iterative scheme, which derived using gradient descent:
\begin{eqnarray}
X^{k+1} = \displaystyle \mbox{argmin}_{X \mbox{~s.t.~} \mathcal{A}(X)=b}~~\mbox{Trace}((X^k + \delta I_N)^{-1}X),
\end{eqnarray}
where $\delta>0$ is a regularization parameter to ensure invertibility of $(X^k + \delta I_N)$, and $X^0$ is typically set to the identity matrix. However, note that this matrix might become il-conditioned with the iterations thereby making the inverse error prone. Moreover, inversion is also computationally intensive as the dimension of the problem grows. \\\\
A variant of \eqref{min_rank_prob} problem is the \textit{General Minimum-Rank-Problem} (MRP) given by:
\begin{eqnarray}
\mbox{Find}~~ &X \in \mathcal{R}^{M\times N} \mbox{~~with ~the ~minimum~ rank} \nonumber \\
\mbox{such~that}~~ &\mathcal{A}(X) = b.
\label{min_rank_prob_gen}
\end{eqnarray}
In the seminal paper \cite{recht2010guaranteed}, it has been proved that minimizing nuclear norm of $X$ results in the minimum rank solution provided the linear operator $\mathcal{A}$ satisfies the Restricted Isometry Property (RIP). The nuclear norm minimization is the following convex program:
\begin{equation}
\begin{array}{l}
\displaystyle \min_{X\in \mathcal{R}^{N\times M}} ~~~\sigma_*(X)\\
\mbox{subject to}\\
\hspace*{1.5cm} \mathcal{A}(X) = b,
\end{array}
\end{equation}
where $\sigma_*(X)$ is the sum of the singular values of $X$. However, checking RIP for a given linear operator is itself NP-Hard in general. In certain special cases where $\mathcal{A}$ is  sampled from Gaussian ensembles, it has been shown that it satisfies RIP with a high probability \cite{recht2010guaranteed}. Loosely speaking, this would mean that given an optimization problem where it is known that the linear constraint operator is sampled from special ensembles, one can find the minimum rank solution with a high probability. In addition, another common method employed to solve  \eqref{min_rank_prob_gen} is the coordinate descent method, in which firstly, $X$ is parameterized as $L\times R$, where $L\in \mathcal{R}^{M\times r}$ and $R\in \mathcal{R}^{M\times r}$. Then $||\mathcal{A}(X) - b||_2$ is minimized iteratively; in each iteration either $L$ or $R$ is alternatively held to the value of the previous iteration. Note that the optimization problem \eqref{min_rank_prob} is devoid of the positive semidefinite constraint. But then one can modify it to the following equivalent problem (as shown in \cite{fazel2003log}):
\begin{eqnarray}
\mbox{Find}~~ &X \in \mathcal{R}^{M\times N}, Y \in \mathcal{R}^{M\times M} \mbox{~and~} Z \in \mathcal{R}^{N\times N} \mbox{~~with ~the ~minimum~  \mbox{Rank}(Y)+\mbox{Rank}(Z)} \nonumber \\
\mbox{such~that}~~ &\mathcal{A}(X) = b 
\mbox{~~and}~~ \begin{bmatrix}
Y&X\\X^\top&Z
\end{bmatrix} \succeq 0.
\label{min_rank_prob_equi}
\end{eqnarray}
The modified problem falls into the class of problems defined by \eqref{main_prob}. Therefore,  \eqref{main_prob} is more general than \eqref{min_rank_prob}, in the sense that if there is an algorithm to solve \eqref{main_prob}, the one can also solve \eqref{min_rank_prob} by iterating over all possible ranks and choosing the one with the lowest rank. Moreover, with this formulation one can also used the \textbf{log-det} heuristic. Also note that minimizing nuclear norm with positive definite constraint boils down to minimizing the trace of the matrix. 

\section{The proposed method}
\begin{figure}[t]
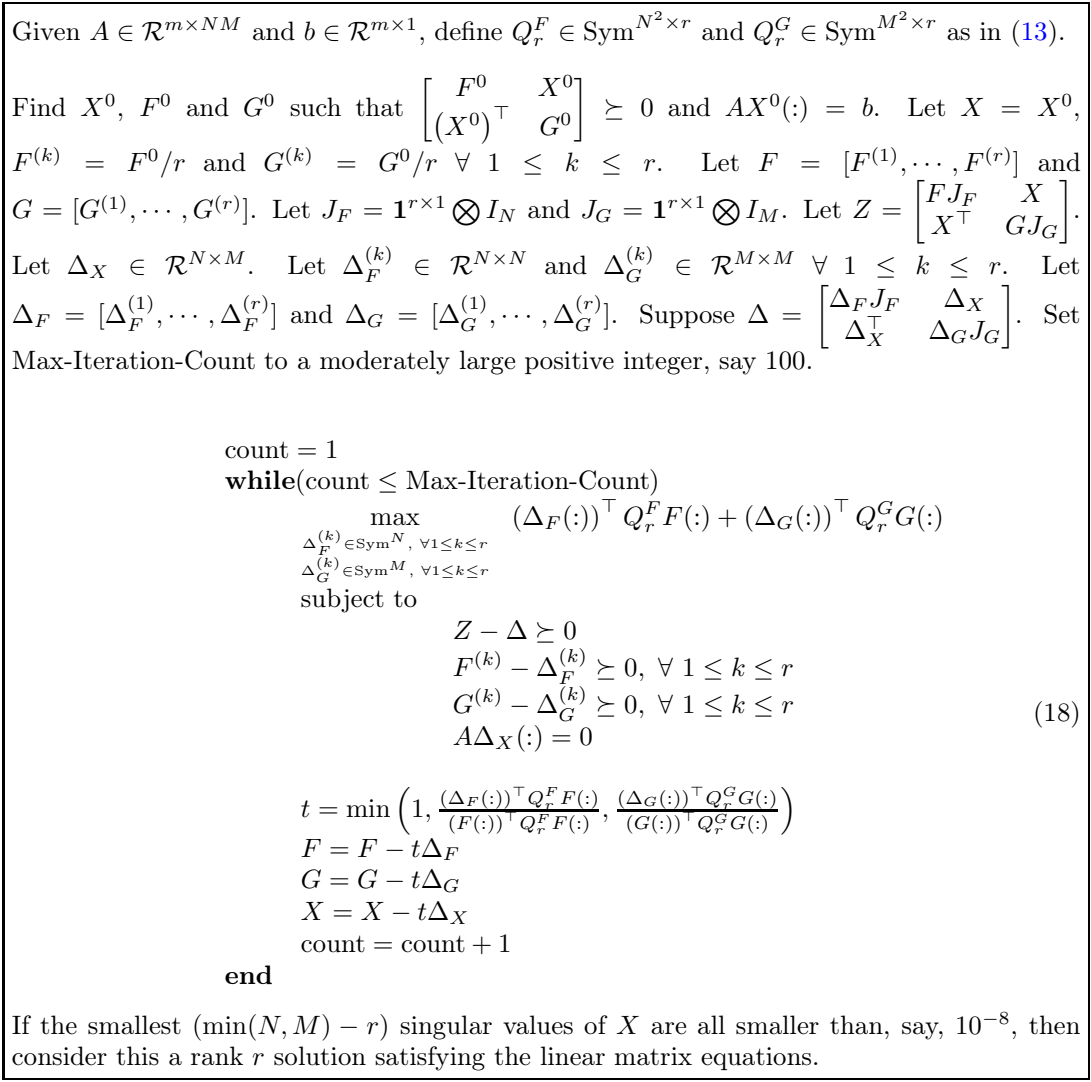

\centering
\fbox{\begin{minipage}{40em}
Given $A \in \mathcal{R}^{m\times NM}$ and $b\in \mathcal{R}^{m\times 1}$, define $Q^F_r\in \mbox{Sym}^{N^2\times r}$ and $Q^G_r\in \mbox{Sym}^{M^2\times r}$ as in \eqref{eqn:part_rank_def}.  \\\\
Find $X^0$, $F^0$ and $G^0$ such that $\begin{bmatrix}F^0 & X^0\\\left(X^0\right)^\top & G^0\end{bmatrix}\succeq 0$ and $AX^0(:)=b$. Let  $X=X^0$, $F^{(k)}=F^0/r$ and $G^{(k)}=G^0/r$ $\forall~ 1\leq k\leq r$. Let $F = [F^{(1)},\cdots,F^{(r)}]$ and $G = [G^{(1)},\cdots,G^{(r)}]$. Let $J_F = \textbf{1}^{r\times 1}\bigotimes I_N$ and $J_G = \textbf{1}^{r\times 1}\bigotimes I_M$. Let $Z = \begin{bmatrix}FJ_F & X\\X^\top & GJ_G\end{bmatrix}$. Let $\Delta_X\in \mathcal{R}^{N\times M}$. Let $\Delta_F^{(k)}\in \mathcal{R}^{N\times N}$ and $\Delta_G^{(k)}\in \mathcal{R}^{M\times M}$ $\forall~ 1\leq k\leq r$. Let $\Delta_F = [\Delta_F^{(1)},\cdots,\Delta_F^{(r)}]$ and $\Delta_G = [\Delta_G^{(1)},\cdots,\Delta_G^{(r)}]$.  Suppose $\Delta = \begin{bmatrix}\Delta_FJ_F & \Delta_X\\\Delta_X^\top & \Delta_GJ_G\end{bmatrix}$. Set Max-Iteration-Count to a moderately large positive integer, say 100.\\\\
\begin{equation}
\begin{array}{l}
\hspace*{1cm}\mbox{count} = 1\\
\hspace*{1cm}\textbf{while} (\mbox{count} \leq \mbox{Max-Iteration-Count}) \\
\hspace*{2cm}\displaystyle \max_{\tiny\substack{\Delta_F^{(k)} \in \mbox{Sym}^N, ~\forall 1\leq k\leq r\\\Delta_G^{(k)} \in \mbox{Sym}^M, ~\forall 1\leq k\leq r}}~~  \left(\Delta_F(:)\right)^\top Q^F_rF(:) + \left(\Delta_G(:)\right)^\top Q^G_rG(:) \\
\hspace*{2cm}\displaystyle \mbox{subject~to}\\
\hspace*{4cm}\displaystyle Z-\Delta \succeq 0\\
\hspace*{4cm}\displaystyle F^{(k)}-\Delta^{(k)}_F \succeq 0, ~\forall~ 1\leq k\leq r\\
\hspace*{4cm}\displaystyle G^{(k)}-\Delta^{(k)}_G \succeq 0, ~\forall~ 1\leq k\leq r\\
\hspace*{4cm}\displaystyle A\Delta_X(:) = 0\\\\
\hspace*{2cm} t = \min\left(1, \frac{\left(\Delta_F(:)\right)^\top Q^F_rF(:)}{\left(F(:)\right)^\top Q^F_rF(:)},\frac{\left(\Delta_G(:)\right)^\top Q^G_rG(:)}{\left(G(:)\right)^\top Q^G_rG(:)}\right)\\
\hspace*{2cm} F = F - t\Delta_F\\
\hspace*{2cm} G = G - t\Delta_G\\
\hspace*{2cm} X = X - t\Delta_X\\
\hspace*{2cm}\mbox{count} = \mbox{count} + 1 \\
\hspace*{1cm}\textbf{end}
\end{array}
\label{rank_one_gradient_general}
\end{equation}
If the smallest $(\min(N,M)-r)$ singular values of $X$ are all smaller than, say, $10^{-8}$, then consider this a rank $r$ solution satisfying the linear matrix equations.
\end{minipage}}
\caption{Heuristic based on gradient descent for finding a rank $r$ constrained solution to \eqref{main_prob_gen}.}
\label{fig:rank_gradient_gen}
\end{figure}

In this section, the proposed heuristic will be outlined. In addition, the intuition leading to the development of the heuristic will also be discussed, a few interesting observations will be made and a few relevant questions will be posed. The author would like to begin by stating the following two standard lemmas.\\\\
\textbf{\textit{Lemma 1}}: For a positive semidefinite matrix $X$, $X_{i,i}X_{j,j}\geq X^2_{i,j}$  for any $1\leq i,j \leq N$.\\
\textbf{\textit{Proof}}: For a positive semidefinite matrix, any principal $2\times 2$ minor is positive definite. Consider the $2\times 2$ principal minor given by $\displaystyle \begin{bmatrix}X_{i,i}& X_{i,j}\\X_{j,i} & X_{j,j}\end{bmatrix}$. Since this is positive semidefinite, its determinant is also positive, that is, $X_{i,i}X_{j,j}-X^2_{i,j}\geq 0$. \hfill \textbf{Q.E.D.}\\\\
\textbf{\textit{Definition 2}}: Note the following definition motivated from \textit{Lemma 1}, termed as the \textit{Low-Rank-Functional} (LRF) $\mathcal{L}: S^N_+\rightarrow \mathcal{R}$: 
\begin{equation}
\displaystyle \mathcal{L}(X) = \sum_{i,j}\left(X_{i,i}X_{j,j}-X^2_{i,j}\right) = X^\top(:)QX(:).
\label{lrf_def}
\end{equation}
A theorem follows immediately from this definition.\\\\
\textbf{\textit{Theorem 3}}: A non-zero positive semidefinite matrix $X$ has unity rank if and only if $\mathcal{L}(X)=0$.\\
\textbf{\textit{Proof}}: Suppose $X$ has unity rank. Then X = $xx^\top$ and therefore $X_{i,j} = x_ix_j$. It is immediately clear that $\mathcal{L}(X) = 0$. The proof of the converse is by induction on the dimension $N$ of the matrix $X$. For $N=2$, it is immediately clear that if for some non-zero matrix $X\succeq 0$,  $\mathcal{L}(X)=det(X)=0$, rank of $X$ has to be unity. Suppose this is true for the case $N\leq M$, for some positive integer $M$ greater than 2. Then one proceeds to prove the statement for dimension $(M+1)$. Suppose $X\in \mbox{Sym}^{(M+1)}$, $X\succeq 0$ and that $\mathcal{L}(X)=0$. Write $\displaystyle X = \begin{bmatrix} X_M & v\\ v^\top & z \end{bmatrix}$. Note that $\displaystyle \mathcal{L}(X) = \mathcal{L}(X_M) + \sum_{i}\left(zX_{i,i}-v_i^2\right)$. Since $\mathcal{L}(X)$ is 0 and  $\displaystyle \sum_{i}\left(zX_{i,i}-v_i^2\right)$ is non-negative (due to Lemma 1), $\mathcal{L}(X_M)$ is also equal to zero. Suppose $X_M$ is a zero matrix. Then, since $\left(zX_{i,i}-v_i^2\right)=0$ for each $i$, $v=0$. Moreover, $z>0$ as it is assumed $X$ is a non-zero matrix. In this case, the rank of $X$ is certainly unity. Now suppose $X_M$ is not a zero matrix. Again, since $\mathcal{L}(X_M)=0$, by the induction hypothesis, $X_M$ has unity rank. Say $X_M = xx^\top$. Now, there are two cases possible: (i) $z = 0$ and (ii) $z>0$. For case (i), $v=0$ and therefore the rank of $X$ still remains unity. For case (ii), by Schur's complement, $X\succeq 0$ if and only if $\displaystyle \left(xx^\top - \frac{1}{z}vv^\top\right) \succeq 0$. This condition implies that $v = kx$, such that $k^2\leq z$. Now suppose $x_r$ is a non-zero element of $x$.  Then, $\left(zx_r^2 - k^2x_r^2\right) = 0$ would imply $z = k^2$. This in turn means that $X = \begin{bmatrix} xx^\top & kx\\kx^\top & k^2 \end{bmatrix}$. Therefore, $X$ has unity rank, and the proof concludes. \hfill \textbf{Q.E.D.}\\\\
\textbf{\textit{Lemma 4}}: For any two positive semidefinite matrices $X$ and $Y$, $\mbox{Tr}(XY) \geq 0$.\\
\textbf{\textit{Proof}}: Consider two unity rank positive semidefinite matrices $A$ and $B$. Since each has unity rank, $A = aa^\top$ and $B = bb^\top$. Now, it is easy to see that $\mbox{Trace}(AB) = (a^\top b)^2 \geq 0$. By Singular Value Decomposition, $\displaystyle X = \sum_{i}^{r_X} X_i$ and $\displaystyle Y = \sum_{i}^{r_Y} Y_i$. Therefore, $\displaystyle \mbox{Trace}(XY) = \sum_{i,j}\mbox{Trace}(X_iY_j)$. Since each of the terms in the summation is non-negative by the last statement on trace of product of the two unity rank PSD matrices, $\mbox{Trace}(XY)\geq 0$. \hfill \textbf{Q.E.D.}\\\\
\textbf{\textit{Definition 5}}: Now for each integer $1\leq r\leq N$, one can make also the following definition, termed here as the \textit{Particular-Rank-Functional} (PRF),  $\mathcal{P}_r:\underbrace{S^n_+\times\cdots\times S^n_+}_{\mbox{r-times}} \rightarrow \mathcal{R}$ as: 
\begin{equation}
\mathcal{P}_r(X^{(1)},\cdots,~X^{(r)}) = \sum_{k,i,j}\left(X^{(k)}_{i,i}X^{(k)}_{j,j}-\left(X^{(k)}_{i,j}\right)^2\right) + \sum_{p\ne q}\mbox{Tr}\left(X^{(p)}X^{(q)}\right) = \tilde{X}^\top Q_r\tilde{X},
\label{eqn:part_rank_def}
\end{equation}
where $\tilde{X}=[X^{(1)}(:);\cdots;X^{(r)}(:)]$. Note that $\mathcal{P}_1(X) = \mathcal{L}(X)$ for a PSD matrix $X$. Now one can derive the following theorem:\\\\
\textbf{\textit{Theorem 6}}: A positive semidefinite matrix $X$ defined as:
\begin{equation}
X = \sum_{k=1}^{r}X^{(k)},
\end{equation}
where each matrix in the summation is also PSD, has rank $r$ if $\mathcal{P}_r(X^{(1)},\cdots,~X^{(r)})=0$. Conversely, if a matrix $X$ has rank $r$, then it can be written as a sum of $r$ positive semidefinite matrices $\{X^{(1)},\cdots,~X^{(r)}\}$, such that $\mathcal{P}(X^{(1)},\cdots,~X^{(r)})=0$.\\
\textbf{\textit{Proof}}: Suppose that for a set of PSD matrices $\{X^{(1)},\cdots,~X^{(r)}\}$, $\mathcal{P}_r(X^{(1)},\cdots,~X^{(r)})$ is zero. Since each of the summands in \eqref{eqn:part_rank_def} is nonnegative (by \textit{Lemma 1} and \textit{Lemma 4}), each term has to be equal to 0. Note that the first summation term in \eqref{eqn:part_rank_def}  is equal to $\displaystyle \sum_{k=1}^{r}\mathcal{L}(X^{(k)})$. Since each of these terms is equal to 0, it  implies that each $X^{(k)}$ has unity rank. Similarly, the summands in the second summation  term in \eqref{eqn:part_rank_def} equal to 0 imply  that the $X^{(i)}$ and $X^{(j)}$ are orthogonal for all $i\neq j$. Hence, the rank of $X = \sum_{k=1}^{r}X^{(k)}$ is equal to $r$. The second statement of the theorem is a corollary of  Singular Value Decomposition. \hfill Q.E.D.\\\\
\textbf{\underline{Construction of the heuristic}}: With the aforementioned results, it is quite intuitive that given an instance of the feasibility problem \eqref{main_prob}, one can replace it by:
\begin{eqnarray} \nonumber
\displaystyle \min_{X^{(1)},\cdots,X^{(k)}}~~& \mathcal{P}_k(X^{(1)},\cdots,X^{(k)})\\
\displaystyle \mbox{such~that}~~ &\mathcal{A}\left(\sum_{j=1}^{k}X^{(j)}\right) = b, ~ X^{(j)} \succeq 0, ~\forall j.
\label{main_prob_heuristic}
\end{eqnarray}
Note that \eqref{main_prob_heuristic} is a non-convex quadratic optimization problem. One standard method to deal with non-convex problems is the gradient descent. This would constitute the first heuristic. Note that the gradient here would be a linear functional which is computationally efficient. The other method is to convert the quadratic problem into a bilinear form and solve a sequence of convex problems fixing one of the variables in the bilinear form to the previous iterate. This would constitute the basic idea for the second heuristic. Note that convex problem in each iteration would be linear semidefinite program which can be solved efficiently using standard optimizers. \\\\
The gradient based heuristic for finding a rank $r$ constrained solution to \eqref{main_prob} is shown in Figure.  \ref{fig:rank_gradient}. For \eqref{main_prob_gen}, a similar gradient based heuristic is shown in Figure. \ref{fig:rank_gradient_gen}. It has been noted that constraining both $F$ and $G$ matrices in \ref{fig:rank_gradient_gen} leads to a rank $r$ solution, although by the semidefinite-embedding theorem (see \cite{fazel2003log}), the rank can be lower than $r$. In both these heuristics, a solution is said to be reached if the smallest $(N-r)$ singular values are lesser than a tolerance, say $10^{-8}$. Similarly, the heuristic based on the bilinear formulation for the problem in \eqref{main_prob} is presented in Figure. \ref{fig:rank_bilinear}. The calculation of SVD before the update step has been observed to improve the convergence rates. Moreover, it has also been observed that if one wants to find a rank $r$ solution and if one starts with an initial guess which has a rank $r+1$, the convergence is quicker. Similarly, the bilinear formulation based heuristic for solving \eqref{main_prob_gen} is presented in Figure. \ref{fig:rank_bilinear_gen}.\\\\
Note that the aforementioned heuristics can also be used with linear matrix inequalities (in this case, $A\mbox{\textbf{vec}}(\Delta J)=0$ is replaced by $A\mbox{\textbf{vec}}(\Delta J) \geq A\mbox{\textbf{vec}}(X J)-b$). Once these  heuristics converge onto a matrix $X$ which has $N-r$ eigenvalues lesser than a tolerance, say $10^{-8}$, then one can perform the following polishing step: compute the SVD of $X = USU^\top$ and then use $U(:,1:r)Z^*$, where $Z^*$ solves
\begin{equation}
\begin{array}{l}
\displaystyle \min_{V \in \mathcal{R}^{r\times N}} ||A\mbox{\textbf{vec}}(U(:,1:r)Z) - b||_2\\
\mbox{subject to}\\
\hspace*{2cm} U(:,1:r)Z \succeq 0.
\end{array}
\label{polish_step}
\end{equation}
It has also been observed that the eigenvectors of the matrix $Q$ in \eqref{lrf_def}, when reshaped into a square matrix are either symmetric or skew-symmetric. For $Q_r$, suppose it that the eigenvectors are divided into $r$ consecutive vectors. Then, it has been observed that each of these parts reshapes either into a symmetric matrix or skew-symmetric matrix. These observations does not have a theoretical basis yet and further properties need to be explored. However, for the matrix $Q$ as defined in \eqref{lrf_def}, the following properties can also be shown easily:
\begin{itemize}
\item For any two PSD matrices $X$ and $Y$ of dimension $N\times N$, $X(:)^\top QY(:)\geq 0$.
\item For any two PSD matrices $X$ and $Y$ of dimension $N\times N$, \newline $(X(:)+Y(:))^\top Q(X(:)+Y(:))\geq X(:)^\top QX(:) + Y(:)^\top QY(:)$.
\end{itemize}
\begin{figure}[t]
\centering
\fbox{\begin{minipage}{40em}
Given $A \in \mathcal{R}^{m\times N^2}$ and $b\in \mathcal{R}^{m\times 1}$, define $Q_r\in \mbox{Sym}^{N^2\times r}$ as shown in \eqref{eqn:part_rank_def}.  \\\\
Find a feasible solution $X^0$ such that $X^0\succeq 0$ and $AX_0(:)=b$. Let  $X^{(k)}=\frac{1}{r}X^0$, for all integers $k$ such that $1\leq k \leq r$. Let $J = \textbf{1}^{r\times 1}\bigotimes I_N$. Let $X = [X^{(1)},\cdots,X^{(r)}]$ and $\Delta = [\Delta^{(1)},\cdots,\Delta^{(r)}]$. Set Max-Iteration-Count to a moderately large positive integer, say 100.\\\\
\begin{equation}
\begin{array}{l}
\hspace*{1cm}\mbox{count} = 1\\
\hspace*{1cm}\textbf{while} (\mbox{count} \leq \mbox{Max-Iteration-Count}) \\
\hspace*{2cm}\displaystyle \max_{\tiny \Delta^{(k)} \in \mbox{Sym}^N, ~\forall 1\leq k\leq r}~~~~  \left(\Delta(:)\right)^\top Q_rX(:) \\
\hspace*{2cm}\displaystyle \mbox{subject~to}\\
\hspace*{4cm}\displaystyle X^{(k)}-\Delta^{(k)} \succeq 0, ~\forall~ 1\leq k\leq r \\
\hspace*{4cm}\displaystyle A\mbox{\textbf{vec}}\left(\Delta J\right) = 0\\\\
\hspace*{2cm} t = \min\left(1, \frac{\left(\Delta(:)\right)^\top Q_rX(:)}{\left(X(:)\right)^\top Q_rX(:)}\right)\\
\hspace*{2cm} X = X - t\Delta\\
\hspace*{2cm}\mbox{count} = \mbox{count} + 1 \\
\hspace*{1cm}\textbf{end}
\end{array}
\label{rank_one_gradient}
\end{equation}
If the smallest $(N-r)$ eigenvalues of $X$ are all smaller than, say, $10^{-8}$, then consider this a rank $r$ solution satisfying the linear matrix equations.   
\end{minipage}}
\caption{Heuristic based on gradient descent for finding a rank $r$ solution to \eqref{main_prob}.}
\label{fig:rank_gradient}
\end{figure}
\begin{figure}[t]
\centering
\fbox{\begin{minipage}{40em}
Given $A \in \mathcal{R}^{m\times N^2}$ and $b\in \mathcal{R}^{m\times 1}$, define $Q_r\in \mbox{Sym}^{N^2}$ as shown in \eqref{eqn:part_rank_def}.  \\\\
Find a feasible solution $X^0$ such that $X^0\succeq 0$ and $AX^0(:)=b$. Set  Max-Iteration-Count to moderately large positive integer, say 100.\\\\
Compute the SVD of $X^0$, \textit{i.e.}, $X^0 = USU^\top$. Define $X^{(k)} = U(:,k)S_{k,k}U^\top(:,k)$ and let $X = [X^{(1)},\cdots,X^{(r)}]$. Let $J = \textbf{1}^{r\times 1}\bigotimes I_N$.\\\\
\begin{equation}
\begin{array}{l}
\hspace*{1cm}\mbox{count} = 1\\
\hspace*{1cm}\textbf{while} (\mbox{count} \leq \mbox{Max-Iteration-Count}) \\
\hspace*{2cm}\displaystyle \min_{Y^{(1)}, \cdots,~Y^{(r)} \in \mathcal{R}^{N\times N}}~~~~ Y(:)^\top Q_r X(:)\\
\hspace*{2cm}\displaystyle \mbox{subject~to}\\
\hspace*{3cm}\displaystyle Y = [Y^{(1)},\cdots,Y^{(r)}]\\
\hspace*{3cm}Y^{(k)}\succeq 0;~~1\leq k\leq r\\
\hspace*{3cm}\displaystyle A\mbox{\textbf{vec}}(YJ) = b\\\\ 
\hspace*{2cm}\mbox{Compute the SVD of}~YJ\mbox{~, that is,~} YJ = USU^\top\\
\hspace*{2cm}X^{(k)} = U(:,k)S_{k,k}U^\top(:,k)\\
\hspace*{2cm}\mbox{count} = \mbox{count} + 1 \\
\hspace*{1cm}\textbf{end}
\end{array}
\label{all_rank_hill_climbing}
\end{equation}
If the smallest $(N-r)$ eigenvalues of $XJ$ are all smaller than, say, $10^{-8}$, then consider this a rank $r$ solution satisfying the linear matrix equations.
\end{minipage}}
\caption{Heuristic based on a bilinear formulation for finding a rank $r$ constrained solution to \eqref{main_prob}.}
\label{fig:rank_bilinear}
\end{figure}
\begin{figure}[t]
\centering
\fbox{\begin{minipage}{40em}
Given $A \in \mathcal{R}^{m\times NM}$ and $b\in \mathcal{R}^{m\times 1}$, define $Q^F_r\in \mbox{Sym}^{N^2}$ and $Q^G_r\in \mbox{Sym}^{M^2}$ as in \eqref{eqn:part_rank_def}. \\\\
Find $X^0$, $F^0$ and $G^0$ such that $\begin{bmatrix}F^0 & X^0\\\left(X^0\right)^\top & G^0\end{bmatrix}\succeq 0$ and $AX^0(:)=b$. Set  Max-Iteration-Count to moderately large positive integer, say 100.\\\\
Compute the SVD of $\begin{bmatrix}F^0 & X^0\\\left(X^0\right)^\top & G^0\end{bmatrix}$, that is, $\begin{bmatrix}F^0 & X^0\\\left(X^0\right)^\top & G^0\end{bmatrix} = USU^\top$. Define $F^{(k)} = U(1:N,k)S_{k,k}U^\top(1:N,k)$ and $G^{(k)} = U((N+1):(N+M),k)S_{k,k}U^\top((N+1):(N+M),k)$. Let $F = [F^{(1)},\cdots,F^{(r)}]$ and $G = [G^{(1)},\cdots,G^{(r)}]$. Let $J_Y = \textbf{1}^{r\times 1}\bigotimes I_N$ and $J_Z = \textbf{1}^{r\times 1}\bigotimes I_M$.\\\\
\begin{equation}
\begin{array}{l}
\hspace*{1cm}\mbox{count} = 1\\
\hspace*{1cm}\textbf{while} (\mbox{count} \leq \mbox{Max-Iteration-Count}) \\
\hspace*{2cm}\displaystyle \min_{\substack{X\in \mathcal{R}^{N\times M},\\ Y^{(1)}, \cdots,~Y^{(r)} \in \mathcal{R}^{N\times N},\\Z^{(1)}, \cdots,~Z^{(r)} \in \mathcal{R}^{M\times M}}}~~ F(:)^\top Q^F_rY(:) + G(:)^\top Q^F_rZ(:)\\
\hspace*{2cm}\displaystyle \mbox{subject~to}\\
\hspace*{4cm}\displaystyle Y = [Y^{(1)},\cdots,Y^{(r)}],~~Z = [Z^{(1)},\cdots,Z^{(r)}]\\
\hspace*{4cm}Y^{(k)}\succeq 0,~~ Z^{(k)}\succeq 0;~\forall~1\leq k\leq r\\
\hspace*{4cm}\displaystyle AX(:) = b\\
\hspace*{4cm}\begin{bmatrix}YJ_Y & X\\ X^\top & ZJ_Z\end{bmatrix} \succeq 0\\
\hspace*{2cm}\mbox{Compute the SVD of}~\begin{bmatrix}F & X\\ X^\top & G\end{bmatrix}\mbox{~, that is,~} \begin{bmatrix}F & X\\X^\top & G\end{bmatrix} = USU^\top\\
\hspace*{2cm}F^{(k)} = U(1:N,k)S_{k,k}U^\top(1:N,k),~\forall~1\leq k\leq r\\
\hspace*{2cm} G^{(k)} = U((N+1):(N+M),k)S_{k,k}U^\top((N+1):(N+M),k), ~\forall~1\leq k\leq r\\
\hspace*{2cm}\mbox{count} = \mbox{count} + 1 \\
\hspace*{1cm}\textbf{end}\\\\
\end{array}
\label{all_rank_general_hill_climbing}
\end{equation}
If the smallest $(\min(N,M)-r)$ singular values of $X$ are all smaller than, say, $10^{-8}$, then consider this a rank $r$ solution satisfying the linear matrix equations.
\end{minipage}}
\caption{Heuristic based on a bilinear formulation for finding a rank $r$ constrained solution to \eqref{main_prob_gen}.}
\label{fig:rank_bilinear_gen}
\end{figure}

\section{Simulations}
In this section, the following specific cases are considered and simulation examples are presented explicitly, that is, (i) finding rank constrained solution to problems of the form given in \eqref{main_prob} and \eqref{main_prob_gen}, and a comparison with the log-det heurisitic and nuclear norm heuristic, (ii) application to Knapsack decision problem, (iii) solving a class of non-linear equations and (iv) application to Fourier Phase Retrieval with amplitude constraints. In additions, two cases where the heuristic fails have also been discussed. Having said this, these cases certainly do not span the entire range of problems which can be tackled with this approach, neither do the cases presented here promise the effectiveness of the proposed method in all scenarios. The author would urge the interested reader to explore the strengths and limitations of the proposed methods. 
\subsection{Rank $r$ constrained solutions}
Here, the problem of finding rank $r$ solutions to LME is considered, with and without the positive semidefinite constraint. For these cases, the heuristics based on the bilinear formulations are used. Note that the heuristic is followed by the polishing step given in \eqref{polish_step}. Two simulation results have been given in the files "PSD\_sim.txt" and "GEN\_sim.txt" uploaded with this document. In both simulations, the matrix $A$, the matrix which was used to generate $b$ (just to ensure that $AX(:)=b$ has at least one solution), the vector $b$ are given at the beginning. In "PSD\_sim.txt" this is followed by the initial guess to start the bi-linear heuristic (here $=0$ implies a zero matrix), results for every rank (and the solution obtained after polishing step), the singular values of the solutions and the norm of $AX(:)=b$. Finally, the results obtained using the trace minimzation and log-det heuristic are also shown. The content of the "GEN\_sim.txt" are similar, except for the fact that the diagonal elements are constrained to 100 here, and that the initial guess for finding rank $r$ solution is the rank $r+1$ solution, if it is found using the proposed heuristic. \\\\
One can note that in "PSD\_sim.txt", solution corresponding to all ranks are obtained using the proposed heuristic, where as the solution obtained using the trace minimization method and log-det heuristic have a rank of 3. On the other hand, one can see that "GEN\_sim.txt" shows solutions corresponding to every rank except unity rank are obtained using the proposed heuristic, where as the nuclear norm minimization leads to a full rank solution. 

\subsection{Knapsack Problem}
The classical problem of Knapsack asks the following question: given a set of objects, each with a weight and a value, what is the subset which yields the maximum sum total value subject to an upper bound on the total weight. The Knapsack decision problem is a minor modification: given a lower bound on the sum total value and an upper bound on the weight, does there exist a subset which satisfies both the constraints. As mentioned earlier, in terms of a rank constrained problem, it can be posed as:
\begin{equation}
\begin{array}{l}
\mbox{Find~} X \in \mathcal{R}^{(N+1)\times (N+1)} \mbox{~with unity rank~}\\
\mbox{subject to~}\\
\hspace*{1cm} X_{1,1} = 1, ~~X \succeq 0,\\
\hspace*{1cm} \mbox{diag}(X) = X_{(:,1)},\\
\hspace*{1cm} \textbf{v}^\top X_{(2:(N+1),1)} \geq V,\\
\hspace*{1cm} \textbf{w}^\top X_{(2:(N+1),1)} \leq W,
\end{array}
\end{equation}
where $\textbf{v}$ and $\textbf{w}$ are the vectors of values and weights respectively, whereas $V$ and $W$ denote the lower bound on total value and upper bound on total weight, respectively. A simulation instance for this is provided.\\\\
    A simulation instance has been provided in the file "Knapsack.txt" provided with this document. In the beginning, it shows that weigths and the values vectors, the minimum total value desired and the maximum allowed total weight and the initial guess for the bilinear heuristic to be applied to this problem. This is followed by the matrix obtained using the proposed heuristic, its singular values, the result for the knapsack problem and the total weight and value of the objects chosen. In addition, it can be seen that the matrix obtained using the trace heuristic has a rank of 15. However, the solution obtained using the log-det heuristic has unity rank, and thus this qualifies as a valid solution.

\subsection{Non-linear Equations}
In this subsection, a solution to a class of non-linear equations (mentioned earlier) is obtained using the proposed method. The problem considered here is:
\begin{eqnarray}
\mbox{Find}~~ &[\theta_1,~\theta_2,\cdots,~\theta_N] \in [-\pi,\pi]^{N\times 1}\nonumber \\
\mbox{such~that}~~ & Ax = b, \mbox{~~where~} x = \left[e^{\mbox{j}\theta_1},~e^{\mbox{j}\theta_2},\cdots,~e^{\mbox{j}\theta_N}\right]^\top.
\label{nl_prob}
\end{eqnarray}
This can be converted to a rank constrained feasibility problem given by:
\begin{eqnarray}
\mbox{Find}~~ &X \in \mathcal{R}^{N\times N}\nonumber \\
\mbox{such~that}~~ &\mathcal{A}X(\small{2:N,1}) = b, ~\mbox{diag}(X) = \textbf{1}, ~X \succeq 0 \mbox{~~and~~} \mbox{Rank} (X) = 1.
\label{combi_prob_conversion}
\end{eqnarray}
A simulation instance for this case is provided in the file "Nonlinear\_Equations.txt" available with this document. In the beginning, it comprises of the matrix $A$, a vector of unit magnitude complex numbers used to generate $b$ (just to ensure that $AX(2:N,1)=b$ has at least one solution) and the initial guess for the proposed bilinear heuristic. This is followed by the solution obtained using the proposed heuristic, its singular values and the solution. Finally, it can be seen that the trace minimization leads to a a full rank solution and log-det heuristic leads to a rank 2 solution.

\subsection{Fourier Phase Retrieval with Amplitude Constraints}
Fourier Phase Retrieval is a classical problem in signal processing, which asks the following: given the magnitude spectrum of a signal (discrete time/frequency), can one find the time domain signal? As the problem as such is ill-posed (has infinitely many solutions, each corresponding to a phase spectrum one chooses), one looks for conditions which renders the solution unique. Here, amplitude constraint are chosen. These constraints need not necessarily lead to a unique solution, but this example serves as a good application for the proposed method. Mathematically, the Fourier Phase Retrieval problem with amplitude constraints can be cast as the following:
\begin{equation}
\begin{array}{l}
\mbox{Find~} X \in \mathcal{R}^{(N+1)\times (N+1)} \mbox{~with unity rank~}\\
\mbox{subject to~}\\
\hspace*{1cm} X_{1,1} = 1, ~~X \succeq 0,\\
\hspace*{1cm} \mbox{diag}(DX(2:(N+1),2:(N+1))D^{\mbox{H}}) = z,\\
\hspace*{1cm} \mbox{diag}(X(2:(N+1),2:(N+1))) \leq A, 
\end{array}
\end{equation}
where $D$ is the DFT matrix of dimension $N\times N$, $z$ is the vector of the magnitudes at different frequencies and $A$ is the amplitude constraint on the time domain signal.\\\\
The simulation instance for this problem has been provided in the file "Phase\_Retrieval.txt" available with this document. It shows the given magnitude squared spectrum (a vector) in the beginning, followed by the amplitude constraint desired. The initial guess used for the bilinear heuristic is shown next. Then the matrix obtained using the proposed heuristic is given, its singular values are shown and the final signal adhering to the magnitude and amplitude constraints is provided. The matrix obtained using trace minimization has full rank and that obtained using log-det has a rank of 4. Hence, both of these do not qualify as a valid solution the Fourier Phase Retrieval problem instance.

\subsection{Shortcomings of the proposed method}
Having mentioned the applications where the proposed method seems to work satisfactorily, it has also been observed that there are applications where the method seems to fail to find a unity rank solution. Two such cases are: (i) the subset sum problem and (ii) the linear complementarity problem in the general symmetric case. However, in both these cases, it has also been observed that the heuristic converges to a rank 2 solution. 
\begin{itemize}
\item {The subset sum problem asks the following question: given a finite set of integers, does there exist a subset whose elements add up to a desired value? Mathematically, the problem can be posed as:
\begin{equation}
\begin{array}{l}
\mbox{Find~} X \in \mathcal{R}^{(N+1)\times (N+1)} \mbox{~with unity rank~}\\
\mbox{subject to~}\\
\hspace*{1cm} X_{1,1} = 1, ~~X \succeq 0,\\
\hspace*{1cm} \mbox{diag}(X) = X(:,1),\\
\hspace*{1cm} \textbf{s}^\top X(2:(N+1),1) = D, 
\end{array}
\end{equation}
where $\textbf{s}$ is the vector of values of the elements in the finite set and $D$ is the desired value.
}

\item{The linear complementarity problem in the general symmetric case is the following problem: given a matrix $M\in \mbox{Sym}^N$ and $q\in \mathcal{R}^{N\times 1}$, find non-negative vectors $w\in \mathcal{R}^{N\times 1}$ and $z\in \mathcal{R}^{N\times 1}$ such that $w^\top z = 0$ and $w = Mz + q$. Mathematically, it can be posed as the following feasibility problem:
\begin{equation}
\begin{array}{l}
\mbox{Find~} X \in \mathcal{R}^{(2N+1)\times (2N+1)} \mbox{~with unity rank~}\\
\mbox{subject to~}\\
\hspace*{1cm} X_{1,1} = 1, ~~X \succeq 0,\\
\hspace*{1cm} \textbf{1}^\top \mbox{diag}(X_{(2:(N+1),(N+2):(2N+1))}) = 0,\\
\hspace*{1cm} X_{(2:(N+1),1)} = MX_{((N+2):(2N+1),1)} + q.
\end{array}
\end{equation}
}
\end{itemize}

\section{Conclusions and Acknowledgements}
The author believes that the method outlined in this paper will pave the way for further research in the domain of rank constrained optimization. The author would also like to extend an invitation to all the interested readers to enrich the proposed method with theoretical guarantees. \\\\
The author wishes to thank Dr. Bharath Bhikkaji of \textit{IIT Madras}, Dr. Arun Ayyar of \textit{Santa Fe Research}, Dr. Girish Ganesan and Dr. Kumar Vijay Mishra for all their support.      

\bibliographystyle{alpha}
\bibliography{sample}

\end{document}